\documentclass[11pt]{article}
\usepackage[letterpaper,width=135mm,height=203mm]{geometry}

\usepackage{graphicx}

\newtheorem{theorem}{Theorem}
\newtheorem{lemma}[theorem]{Lemma}
\newtheorem{question}{Question}
\newenvironment{Proof}{\textsc{Proof.}}{\hspace{8mm}\QEDmark\smallskip}
\newcommand{\QEDmark}{\mbox{\textsc{qed}}}

\newcommand{\dd}[1]{\textbf{\textit{#1}}}

\newcommand{\curlyP}{\mathcal{P}}
\newcommand{\curlyD}{\mathcal{D}}
\newcommand{\curlyE}{\mathcal{E}}
\newcommand{\curlyF}{\mathcal{F}}

\newcounter{myfig}
\newenvironment{myfigure}{\refstepcounter{myfig}\begin{center}}{\end{center}}
\renewcommand{\caption}[1]{\textsf{Figure \arabic{myfig}:} #1}

\begin{document}

\begin{center}
\textbf{\Large The Generalized Matcher Game} \\[4mm]
\textsuperscript{1}Anna Bachstein, \textsuperscript{1}Wayne Goddard, \textsuperscript{2}Connor Lehmacher\\[2mm]
\textsuperscript{1}School of Mathematical and Statistical Sciences, Clemson University\\[2mm]
\textsuperscript{2}Department of Mathematics, Vanderbilt University
\end{center}

\begin{abstract}
Recently the matcher game was introduced. In this game, two players
create a maximal matching by one player repeatedly choosing a vertex and
the other player choosing a $K_2$ containing that vertex. One
player tries to minimize the result and the other to maximize the result.
In this paper we propose a generalization of this game where $K_2$ is replaced
by a general graph $F$. We focus here on the case of $F=P_3$.
We provide some general results and lower bounds
for the game, investigate
the graphs where the game ends with all vertices taken, and calculate
the value for some specific families of graphs. 
\end{abstract}

\section{Introduction}

Fix some graph $F$ with a designated ``root'' vertex $r$.
Given a graph $G$, two players take turns.
One player \dd{initiates} by choosing a vertex $v$, subject to the constraint
that $G$ contains at least one (not necessarily induced) copy of
$F$ with vertex $v$ corresponding to vertex $r$. The other player 
\dd{responds} by choosing one such copy of $F$ within $G$.
Vertices can only be used once. 
This process continues until the remaining vertices of $G$ do not contain a 
copy of $F$.

One player tries to maximize the number of copies taken.
The other player tries to minimize this number. We call these players
\dd{Maximizer} and \dd{Minimizer}. Thus there are
two versions, depending on who initiates and who
responds. We define the \dd{value} of the game as the number of copies
taken with optimal play by both players.

This game is a generalization of the \dd{matcher game} introduced in \cite{GH}.
That game is where  the subgraph $F$ is $K_2$.
For example, it was shown in \cite{GH} that if Maximizer is the responder,
then the value of the game is just the matching number of the graph.

We focus on the graph $P_3$. Rooted at the center vertex we call it
the \dd{2-star} or simply the \dd{star}; rooted at an end-vertex we call it the \dd{stripe}.
Thus we talk of the \dd{star-game} and the \dd{stripe-game}.

We proceed as follows. In Section 2 we provide some examples and elementary results.
In Section 3 we determine lower bounds
for the game in general graphs and in Section 4 we
consider the graphs where the game ends with all vertices taken.
In Section 5 we consider some specific families of graphs including
grids. Finally in Section 5 we consider the alternative game where there 
is no ``root'' vertex. 

\section{Examples}

For a graph $G$, we define a \dd{$P_3$-packing} as a collection of vertex-disjoint 
copies of $P_3$ in $G$. Further, we denote by $\mu(G)$ the maximum size
of a $P_3$-packing of~$G$.
This parameter generalizes the matching number and is well-studied.
For example, Kaneko et al.~\cite{KKN} showed that if $G$ is a $2$-connected
claw-free graph of order $n$ a multiple of $3$, then $\mu(G) = n/3$.
Earlier, Kirkpatrick and Hell~\cite{KH} showed that the parameter is NP-complete to compute.

The parameter $\mu(G)$ provides an immediate upper bound on the value of
the star- or stripe-game. At the other extreme is the minimum size of a maximal
$P_3$-packing. It is easy to see that that quantity is at least $\mu(G)/3$. 

As a first example, consider the complete bipartite graph. Note that 
up to symmetry, the response is forced. So Maximizer as initiator
can ensure a maximum $P_3$-packing and Minimizer as initiator can ensure
a minimum maximal $P_3$-packing, regardless of whether it is the star- or stripe-game. Thus:

\begin{lemma}
Consider the complete bipartite graph $K_{r,s}$ with $r\le s$ and $s\ge 2$. \\
The value of the game with Maximizer initiating is $\mu = \min( r, \lfloor (s+r)/3 \rfloor )$. \\
The value of the game with Minimizer initiating is $\lceil r/2 \rceil$.
\end{lemma}

Consider next the game played on a path. If Maximizer initiates, they
ensure (almost) all the vertices by choosing an end-vertex for the 
stripe-game and a neighbor of an end-vertex for the star-game.
So we consider the version where Maximizer responds.

For the star-game, there is a unique $P_3$ for a given central vertex.
Thus the game where Minimizer initiates is equivalent to minimum maximal $P_3$-packing.
The arrangement is to skip two vertices, take a $P_3$, skip two vertices, etc.
Thus the value of the star-game played on $P_n$ is $\lfloor (n+3)/5 \rfloor$ if Minimizer initiates.  

\begin{lemma}
Consider the stripe-game with Maximizer responding.
The value of the game played on the path $P_n$ is $\lfloor (n+1)/4 \rfloor$.
\end{lemma}
\begin{Proof}
Assume the vertices are numbered from $1$ up to $n$.
Minimizer can ensure at least the claimed value by playing in succession vertex
number $2$, $6$, $10$, etc. If $n \equiv 3$ mod 4, then a final initiation of vertex $n-1$ 
is invalid, but (both) $n$ and $n-2$ are valid final initiations. 

To show that Maximizer as responder can ensure at most the claimed value,
it suffices to show by induction that the recurrence relation
\[
    f(n) = \min_{1\le k\le n}  \max\{ f(k-1)+f(n-k-2),\, f(k-3)+f(n-k) \} ,
\]
has solution $g(n) = \lfloor (n+1)/4 \rfloor$. For $k'=k+4$ it holds that
$g(k'-1)+g(n-k'-2) = g(k-1)+g(n-k-2)$ and that 
$g(k'-3)+g(n-k') = g(k-1)+g(n-k)$.
So it suffices to check the recurrence for say $1\le k \le 4$. 
From this it follows that if the recurrence is true for $n$ then it is 
also true for $n+4$. So it suffices to check the recurrence for four consecutive
values of $n$, e.g.{} $4\le n \le 7$. This can be performed by hand or computer.
\end{Proof}

\section{General Lower Bounds}

\subsection{Lower bounds for Maximizer responding}

\begin{lemma} \label{l:starMaxRes}
Consider the star-game with Maximizer responding.
If $G$ is a graph with $\mu(G)=m$, then
the value of the game is at least $\lceil m/2\rceil$ and this is sharp. 
\end{lemma}
\begin{Proof}
For bound: consider a maximum $P_3$-packing $\curlyP$ of $G$.
Maximizer as responder can ensure that each move overlaps at most
two of these copies. If the initiation vertex $v$ is outside $\curlyP$, this is
immediate; if $v$ is chosen in some copy of~$\curlyP$, then there is an edge to a neighbor
of $v$ within the copy, and Maximizer can use that neighbor. Thus the game
lasts at least $\lceil m/2\rceil$ moves.

For optimality:  take $m$ copies of $P_3$, and
pick one end-vertex from each copy and make all the chosen
vertices into a path $K$. See Figure~\ref{f:comb} where the vertices of~$K$ are drawn in white.
Minimizer initiates at a vertex in~$K$;
Maximizer responds with a star using another vertex of $K$ 
and one degree-$2$ vertex.
\end{Proof}

\begin{myfigure}
\includegraphics{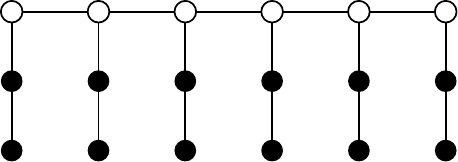}\\
\caption{A $P_3$-packable graph with smallest value of game for Maximizer responding}
\label{f:comb}
\end{myfigure}

\begin{lemma}
Consider the stripe-game with Maximizer responding.
If $G$ is a graph with $\mu(G)=m$, then
the value of the game is at least $\lceil m/2\rceil$ and this is sharp. 
\end{lemma}
\begin{Proof} 
For bound: consider a maximum $P_3$-packing $\curlyP$ of $G$.
We claim Maximizer as responder can ensure that each move overlaps at most
two copies in~$\curlyP$. Suppose the initiation is the end-vertex $a$ of a stripe $abc$ that
intersects three copies $Q_a$, $Q_b$, $Q_c$ of $P_3$ in $\curlyP$. Then vertex $b$ 
has a neighbor $d$ within $Q_b$, so Maximizer can choose the stripe $abd$ instead.
Again the game lasts at least $\lceil m/2\rceil$ moves.

For optimality: use the same construction as in Figure~\ref{f:comb} above.
Minimizer plays a vertex adjacent to a leaf;
Maximizer is forced to respond with a stripe using two white vertices.
\end{Proof}

\subsection{Lower bounds for Minimizer responding}

\begin{lemma} \label{l:starMinRes}
Consider the star-game with Minimizer responding.
If $G$ is a graph with $\mu(G)=m$, then
 the value of the game is at least $\lceil m/3\rceil$ and this is sharp. 
\end{lemma}
\begin{Proof}
For bound: consider a $P_3$-packing of $G$ with $m$ copies.
No matter what the players do, each move can overlap at most three of
these copies. So the game lasts at least $ m/3$ moves. In other words,
the lower bound follows from the lower bound on the size of a
maximal $P_3$-packing noted earlier. 

For optimality: consider the ``double corona'' of the complete graph $K_m$,
with $m$ a multiple of $3$.
That is, take $m$ copies of $P_3$ and
add an edge between every two center vertices to form clique $K$. See Figure~\ref{f:doubleCorona}.
The initiator has to choose a vertex from $K$. Minimizer can respond by taking two more 
vertices from $K$. The game stops after $\lceil m/3\rceil$ moves.
\end{Proof}

\begin{myfigure}
\includegraphics{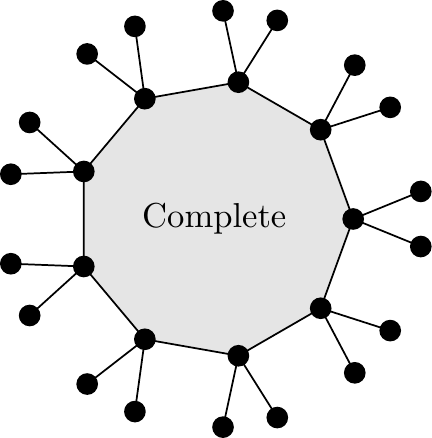}\\
\caption{A $P_3$-packable graph with smallest value of star-game for Minimizer responding}
\label{f:doubleCorona}
\end{myfigure}

The above result can be improved for a tree:

\begin{lemma}  \label{l:starMinRespTree}
Consider the star-game with Minimizer responding.
If $T$ is a tree with $\mu(T)=m$, 
then the value of the game is at least $\lceil m/2\rceil$ and this is sharp. 
\end{lemma}
\begin{Proof}
For bound: consider an optimal $P_3$-packing $\curlyP$ of $G$.
Consider a vertex~$v$ that has at most one non-leaf neighbor. 
Maximizer initiates at vertex $v$. Then $v$ has at most one edge
to a non-leaf, and so Minimizer can overlap at most two
stars in $\curlyP$ when responding. Thus the game lasts at least $m/2$ moves.

Optimality is achieved by the caterpillar formed from $m$
$2$-stars by adding edges so that their centers form a path. See Figure~\ref{f:caterpillar}.
\end{Proof}

\begin{myfigure}
\includegraphics{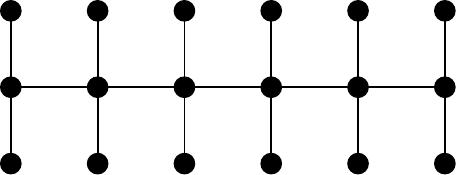}\\
\caption{A $P_3$-packable tree with smallest value of star-game for Minimizer responding}
\label{f:caterpillar}
\end{myfigure}

Consider the stripe-game with Minimizer responding.
If $G$ is a graph with $\mu(G)=m$, then by the same trivial
argument as before the  value of the game is at least $\lceil m/3\rceil$. 
It seems unlikely that this bound is achievable in general, but we
are unable to prove this. In general, we pose the question:

\begin{question}
What is the best possible lower bound for the value of the star-game with Minimizer
responding as a function of the $P_3$-packing number?
\end{question}

As in the case of the star-game, one can prove a better bound if the graph is a tree:
 
\begin{lemma}
Consider the stripe-game with Minimizer responding.
If $T$ is a tree with $\mu(T)=m$, 
then the value of the game is at least $\lceil m/2\rceil$ and this is sharp. 
\end{lemma}
\begin{Proof}
For bound: consider an optimal $P_3$-packing $\curlyP$ of $T$.
Maximizer initiates on an end-vertex $v$. 
Then when Minimizer
responds, they have only one option for the neighbor of $v$, and so 
the stripe taken can overlap at most two copies from~$\curlyP$.

Optimality: Consider the caterpillar in Figure~\ref{f:caterpillar}.
When $m$ is even,
no matter where Maximizer initiates, Minimizer can respond by using
two vertices of the spine and leaving the remaining vertices of the spine
to induce paths each with an even number of vertices.
The argument when $m$ is odd is similar.
\end{Proof}


\section{Perfect Graphs}

In the previous section we considered graphs where the value of the game is small.
At the other end of the spectrum are those graphs where the number
of copies taken is the largest it can be. Define a graph as \dd{perfect} if the 
result of the game is that all vertices 
are taken. 

If a graph has a $K_3$-packing (meaning one can partition the vertex
set into triples such that each induces a triangle), then it is immediate that Maximizer 
as responder can get all the vertices, whether it is the star- or stripe-game.
For example, it is easy to build cubic graphs that have a $K_3$-packing 
(just
start with a collection of triangles and add a perfect matching).

So we ask for trees: which trees are perfect? In the following we 
use the fact that a disconnected graph is perfect if and only if each component
is perfect.

\subsection{Perfect trees with Maximizer as responder}

Let $\curlyD$ be the family of forests defined as follows.
Take some number of disjoint $2$-stars and then add edges
between their centers without creating a cycle.
(One might call each component the double-corona of a tree.) 
This includes for example the caterpillar shown in Figure~\ref{f:caterpillar}.

\begin{lemma}
Consider the star-game with Maximizer responding. 
Then the trees in $\curlyD$ are the perfect trees.
\end{lemma}
\begin{Proof}
The graphs in $\curlyD$ are perfect, since Minimizer is forced each time to 
initiate on a center vertex $v$, and Maximizer can respond
by taking $v$ and its two leaves, leaving a member of $\curlyD$.

We argue that these are the only perfect trees. Consider a perfect
tree $T$ with an initiation
by Minimizer on vertex $a$ and Maximizer's chosen response using vertices $b$ and $c$.
Suppose that $b$ is a valid initiation vertex for Minimizer,
and Maximizer's response to that would include the new vertex $x$. Let
$C_x$ be the component of $T-\{a,b,c\}$ containing~$x$; this 
is by assumption $P_3$-packable and so has order a multiple of~$3$.
But if Minimizer initiates on $b$ and Maximizer responds using $x$,
then this uses one vertex from $C_x$ and so what remains does
not have order a multiple of $3$, and so the tree is not perfect.
That is, any potential initial vertex $v$ must have two leaf neighbors.
Repeat. (Note that the above argument shows that all of $v$'s
other neighbors must have two leaf neighbors, since they too are
potential initial vertices.) 
\end{Proof}

\begin{lemma}
Consider the stripe-game with Maximizer responding. 
Then $P_3$ itself is the only perfect tree.
\end{lemma}
\begin{Proof}
If a tree $T$ has more
than three vertices, then it contains a non-leaf vertex $v$ that has 
exactly one non-leaf neighbor. 
Minimizer initiates on $v$, so that $v$'s leaf neighbors immediately 
become isolated.
\end{Proof}

\subsection{Perfect trees with Minimizer as responder}

Let $\curlyE$ be the family of forests defined as follows.
Start with some number of $P_3$'s. Repeatedly add a $P_3$ and add at most
one edge between it and each existing component, except no edge is
added incident with the central vertex of the new~$P_3$. A member of 
$\curlyE$ is draw in Figure~\ref{f:nineTree}.

\begin{myfigure}
\includegraphics{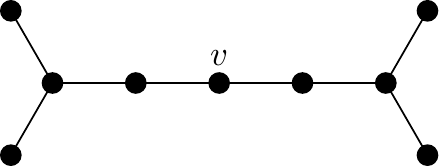}\\
\caption{A perfect tree for star-game with Minimizer responding}
\label{f:nineTree}
\end{myfigure}

\begin{lemma}
Consider the star-game with Minimizer responding. 
Then the trees in $\curlyE$ are the perfect trees.
\end{lemma}
\begin{Proof}
The graphs in $\curlyE$ are perfect,
since Maximizer can initiate on the central vertex $v$ of the final added $P_3$
and then use recursion. (Minimizer is never given a choice.)

We argue that these are the only perfect trees. Consider a perfect
tree $T$ with initiation
by Maximizer on a vertex $a$ and assume one possible response is the star $bac$.
Suppose $a$ has degree more than $2$ 
and let $x$ be one of $a$'s remaining neighbors.
Since the component $C_x$ of $T-\{a,b,c\}$ containing $x$  
is $P_3$-packable, it has order
a multiple of $3$. If instead Minimizer plays $bax$, 
then we have still isolated $C_x$ but removed one vertex from it, so it
does not have order a multiple of~$3$ any more, a contradiction.
Thus we have shown that Maximizer must initiate on a vertex of degree $2$. 
After removal of the star, apply induction.
\end{Proof}

Let $\curlyF$ be the family of forests defined as follows.
Start with nothing. Repeatedly add a $P_3$
with a designated end-vertex $v$, and join $v$ 
to at most one vertex in each existing component. An example
is shown in Figure~\ref{f:curlyF}, where the designated end-vertices
are numbered in order of creation.

\begin{myfigure}
\includegraphics{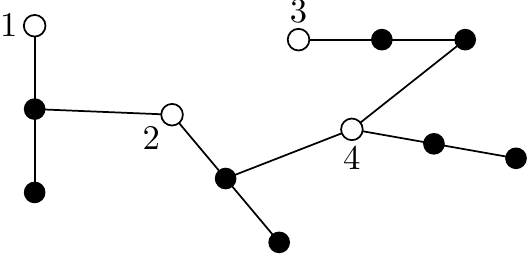}\\
\caption{A perfect tree for stripe-game with Minimizer responding}
\label{f:curlyF}
\end{myfigure}

\begin{lemma}
Consider the stripe-game with Minimizer responding. 
Then the trees in $\curlyF$ are the perfect trees.
\end{lemma}
\begin{Proof}
The graphs in $\curlyF$ are perfect, since 
Maximizer can initiate on the end-vertex of the final added $P_3$ that is not $v$,
and then use induction/recursion. (Minimizer is never given a choice.)

We argue that these are the only perfect trees. Consider a
perfect tree $T$ with initiation by Maximizer on a vertex $a$ and 
assume one possible response by Minimizer is the stripe $abc$.
Suppose $a$ is not an end-vertex
and let $x$ be one of $a$'s remaining neighbors.
Let $C_x$ be the component of $T-\{a,b,c\}$ containing $x$. 
Since it is coverable, $C_x$ has order a multiple of $3$; further 
$x$ has another neighbor, say $y$.
Consider the result if Minimizer responds 
by taking $axy$. This removes two of the vertices from $C_x$; and
so what remains of $C_x$, does not have order a multiple of $3$,
a contradiction. That is, $a$ must be a leaf.

Suppose now that $b$ has 
degree more than $2$; say with another neighbor $z$. Since the component
$D_y$ of $T-\{a,b,c\}$ containing $z$ is $P_3$-coverable,
it has order a multiple of $3$. But if the opening move is $abz$,
then we have still isolated $D_y$ but removed one vertex from it, so it
does not have order a multiple of $3$ any more, a contradiction.

Thus we have shown that Maximizer must initiate on a vertex $a$ 
such that $a$ is end-vertex and its neighbor $b$ has degree $2$.
After removal of $abc$, apply induction.
\end{Proof}

\subsection{Maximal outerplanar graphs}

Recall that a \dd{maximal outerplanar graph} (MOP) is created by taking
a cycle and adding noncrossing chords until their addition
 is impossible. These graphs are a subset of the $2$-trees.
 We consider which such graphs are perfect.

The triangle $K_3$ is always perfect. There are three MOPs of order $6$:
the fan, snake, and sun, drawn here.
\begin{myfigure}
\begin{tabular}{ccc}
\includegraphics{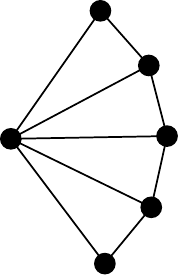} \qquad&\qquad
\includegraphics{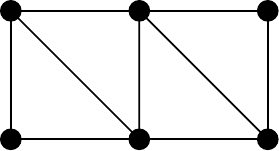} \qquad&\qquad
\includegraphics{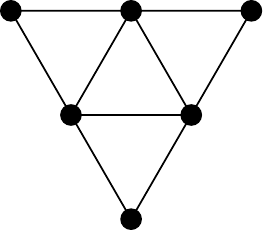}
\end{tabular}\\
\caption{The three MOPs on $6$ vertices}
\end{myfigure}

Perhaps surprisingly, which is perfect does not depend on the graph but only on the game.
Minimizer responding in the stripe-game can ensure the value is~$1$;
In the other three games, each graph is perfect. We omit the details.
The imperfection generalizes:

\begin{lemma}
Consider the stripe-game with Minimizer responding.
Then the only perfect maximal outerplanar graph is $K_3$.
\end{lemma}
\begin{Proof}
Assume the graph is perfect and Maximizer initiates at vertex $v$.
Suppose that $v$ has 
degree more than $2$ and let $vw$ be a chord incident with $v$.
Then removal of $\{v,w\}$ separates the graph $G$ into two components;
let $x$ be a neighbor of~$w$ on the outer cycle, chosen in the component of order a multiple
of $3$ if there exists such a component.
Then Minimizer plays the stripe $vwx$ and leaves neither component a 
multiple of $3$,
and hence not all vertices are eventually taken. That is, the graph is not perfect.

So vertex $v$ has degree $2$. Say its
neighbors are $y_1$ and $y_2$. Then $y_1$ and $y_2$ have another common
neighbor, say $z$. The removal of stripe $vy_1z$ creates a component $C$
without $y_2$;
changing to the removal of stripe $xy_2z$ alters the size of $C$ by 
$1$. So at least one of these removals
produces a component not a multiple of $3$, which gives Minimizer a 
suitable response to avoid all vertices being taken. That is, the graph is not perfect.
\end{Proof}

It is unclear what the perfect MOPs look like for the other three
games. We do note that the first part of the above proof 
carries over to the star-game on a MOP with Minimizer responding: to have a chance
of a perfect outcome, Maximizer must initiate on a vertex of degree $2$.


\section{Some Grid-Like Graphs}

\subsection{Grids with two rows}

\begin{lemma}  \label{l:grid.75}
Consider the star-game. The value of the game on 
a $2\times m$ grid is $\lceil m/2 \rceil$, regardless of which player has which role. 
\end{lemma}
\begin{Proof}
Think of the grid as $2$ rows and $m$ columns.
Consider first the game with Maximizer initiating.

Maximizer as initiator can ensure the desired quantity.
Play top row first column, then top row third column, and so on.
 If $m$ is odd, add one final move of 
bottom row, last column. Each response is forced: each initiated vertex has exactly two neighbors
at the time of being chosen. See Figure~\ref{f:gridOne}.

\begin{myfigure}
\includegraphics{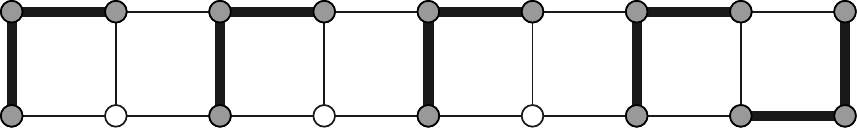}\\
\caption{A maximal packing of size $5$ in the $2\times 9$ grid}
\label{f:gridOne}
\end{myfigure}

Minimizer as responder can ensure at most the desired quantity. 
Assume first that $m$ is even.
Minimizer will always use a vertical edge; further, 
if the star is initiated in column $i$ then: if $i$ is odd, 
they use the vertex in column $i+1$; and if $i$ is even,
they use the vertex in column $i-1$. Equivalently, Minimizer partitions the 
grid into $2\times 2$ subgrids, and responds to an initiation 
in some $2\times 2$ subgrid by staying within that subgrid. 
Note that the fourth vertex of the $2\times 2$ subgrid can never thereafter
be chosen by the
initiator.

If $m$ is odd, Minimizer plays the same strategy where possible.
Specifically,
Minimizer tentatively partitions the grid into $2\times 2$ subgrids with 
one ``floating'' column in the last column. 
If Maximizer initiates in the last column, then Minimizer uses the vertical edge and one vertex of column $m-1$,
as forced. Mentally, Minimizer slides the floating column two to the left and continues the strategy.
If initiator plays in the floating column a second time, then again Minimizer uses the vertical edge and the vertex to the left,
and slides the floating column two to the left. Eventually the floating column will be surrounded
by played $2\times 2$ subgrids. See Figure~\ref{f:gridTwo} for an example.
Thus the number of moves is at most one more than the number 
of $2\times 2$ grids, which equals $\lceil m/2 \rceil$.

\begin{myfigure}
\includegraphics{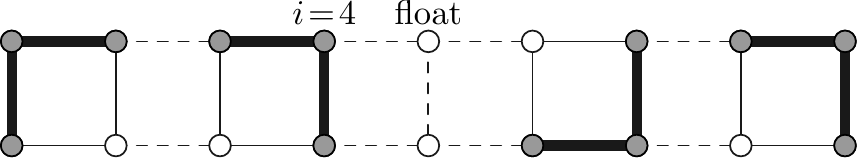}\\
\caption{Play in the $2\times 9$ grid}
\label{f:gridTwo}
\end{myfigure}

The analysis of the game with the roles reversed is the same!
Minimizer initiating can use the same strategy to end the game in $\lceil m/2 \rceil$ moves.
Maximizer responding can use the same strategy to ensure it lasts at least $\lceil m/2 \rceil$ moves.
\end{Proof}

\begin{lemma}
Consider the stripe-game. The value of the game on 
a $2\times m$ grid is $\lceil m/2 \rceil$ if Maximizer responds
and $\lfloor m/2 \rfloor$ if Minimizer responds.
\end{lemma}
\begin{Proof}
Minimizer as initiator can ensure at most  $\lceil m/2 \rceil$ moves as follows. Play
in the top left corner. If Maximizer's stripe is horizontal, then play in the bottom
row in the second column and repeat the strategy as if the first four columns are gone.
If Maximizer's stripe uses a vertical edge, then repeat as if the first two columns are gone.
See FIgure~\ref{f:gridThree}.
Similarly, Maximizer as initiator can ensure at least  $\lfloor m/2 \rfloor$ moves by the same
strategy.

\begin{myfigure}
\includegraphics{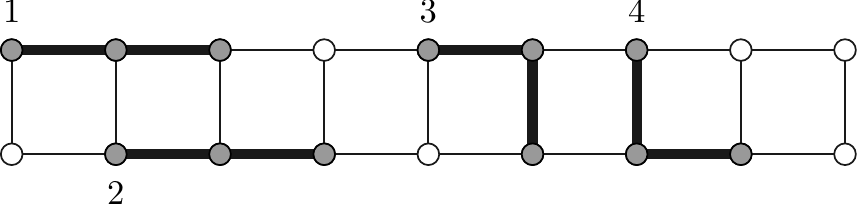}\\
\caption{Start of a stripe-game in the $2\times 9$ grid}
\label{f:gridThree}
\end{myfigure}

Now, we argue that Minimizer as responder
can ensure at most $\lfloor m/2 \rfloor$ moves.
Their strategy is as follows. Number the columns $1$, $2$, up to $m$.
Every response, the Minimizer:
\begin{quote} \slshape
    uses two vertices from an even-numbered column and 
    one vertex from an odd-numbered column. 
 \end{quote}

We claim Minimizer can always achieve this.
If Maximizer initiates with $v$ in an odd-numbered column, 
and neither even-numbered column next to it is available, then 
$v$ is not a valid start vertex.
If Maximizer initiates with $v$ in an even-numbered column, 
and Minimizer cannot respond by starting with the vertical edge, 
then that means the neighboring odd-numbered columns each have 
a vertex taken (or don't exist), and so the even-numbered columns next over are,
by the strategy, taken  (or don't exist), 
and so this is again not a valid move for Maximizer (instead $v$ is the center of a star with 3 leaves). 

There are $\lfloor m/2 \rfloor$ even-numbered columns. And so that is
an upper bound on the number of moves the Minimizer can be forced to make.

Maximizer as responder can ensure at least $\lceil m/2 \rceil$ moves.
This uses a similar idea to above. Specifically:
\begin{quote} \slshape
    every response uses two vertices from an odd-numbered column and 
    one vertex from an even-numbered column. 
 \end{quote}
 As above, if such a response is not possible, then the initiator chose an invalid vertex.
 Further, we claim that if there is no valid vertex available for the initiator, then every odd-numbered
 column is full. The lower bound and the result follows.
 \end{Proof}

 \subsection{Grids with three or more rows}
 
 Maximizer can do well sometimes on the grid with three rows.
  
  \begin{lemma}
For the star-game with Minimizer responding, the three-row grid is perfect.
\end{lemma}
\begin{Proof}
Assume we have a grid with three rows and with columns numbered from $1$ up to $m$.
There are two distinct strategies based on the parity of $m$.

 \emph{Even case:}
Maximizer starts in the lower left corner. 
This gives Minimizer a forced move. Then, Maximizer plays the bottom    
vertex in the third column, forcing the Minimizer again. 
This continues with Maximizer playing the bottom vertex 
in column $2i+1$ for increasing $i$. Thereafter, Maximizer 
plays in the top right corner, again forcing the Minimizer, 
and then moves left across the top row initiating in column $2i+1$ for decreasing $i$.

\emph{Odd case:} Here Maximizer starts with the middle vertex in the first column. If Minimizer responds within the column, then 
we are back in the even case. So 
without loss of generality assume that Minimizer responds by using the vertex in the second column
and the vertex in the bottom left hand corner. Then Maximizer initiates in the top row second column.
Minimizer has a forced response. Maximizer continues by playing in column $2i$ in the top row for increasing
$i$, followed by playing along the bottom row in column $m-2i$ for decreasing $i$. See 
Figure~\ref{f:gridFour}.
Each response after the first move is forced.
\end{Proof}

\begin{myfigure}
\includegraphics{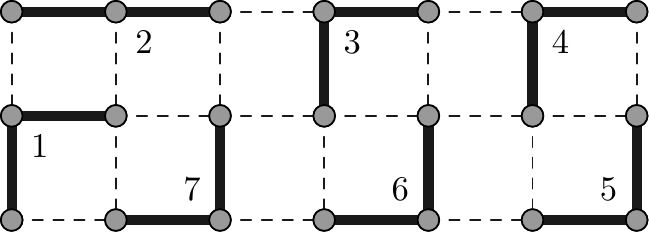}\\
\caption{A perfect star-game}
\label{f:gridFour}
\end{myfigure}
   
 It can be shown that for the stripe-game on three-row grid with $m$ columns, with Maximizer
 as initiator, the value is $m$ if $m\le 3$, and $m-1$ otherwise. In contrast,
 the value of either game when Minimizer initiates is asymptotically at most $(1-\varepsilon)m$
 for some $\varepsilon>0$. We omit the proofs. (In particular we do not know the optimal $\varepsilon$.)

For general grids, we note that initiator can use the same strategy as in Lemma~\ref{l:grid.75} to 
achieve $3/4$ the vertices being used. But it is unclear whether this is good, bad, or indifferent.

\subsection{Rooks graphs with two rows}

We consider \dd{rooks graph} with two rows. That is, the graph $R_m$
obtained by taking two cliques of size $m$ and adding a perfect matching
between them, which we call \dd{cross-edges}. For simplicity, we restrict to $m$ a multiple of $3$.

Every maximal $P_3$-packing for such graph will leave at most $3$ vertices;
and if so, the residue will be $2$ vertices from one clique and $1$ from the other.
Thus the value of the game is either $m$ or $m-1$.

If Maximizer is the responder, then they can stay within the row, 
maintaining each row's size a multiple of $3$, and therefore always get every vertex. That~is:

\begin{lemma}
If Maximizer is responding, then for both games $R_m$ is perfect for all $m$ a multiple of $3$.
\end{lemma}

So we 
consider the game where Minimizer responds.

\begin{theorem}
Consider the star-game with Minimizer responding. 
Then $R_m$ is perfect for all $m$ a multiple of $3$.
\end{theorem}
\begin{Proof}
Maximizer initiates somewhere. (The graph is vertex-transitive.) There are two
cases.

\dd{Case~1:} \textit{Minimizer stays within that row.}
Then, Maximizer initiates at any vertex that does not have
a cross-edge. This forces Minimizer
to play within the row, and ensures that the number of vertices
in each row remains a multiple of~$3$. 
If every vertex has a cross-edge, then we are back 
to a rooks graph, and can apply induction; otherwise Maximizer
continues with a vertex that does not have a cross-edge.

\dd{Case~2:} \textit{Minimizer uses the cross-edge for the first star.} 
Say Maximizer initiated in the top row,
so that two vertices were taken from the top row and one from the bottom row.
Then Maximizer plays the vertex in the bottom row that has no cross-edge,
to which Minimizer is forced to respond by taking three in the bottom row.
Thereafter, Maximizer plays a vertex in the bottom row,
and repeats so long as Minimizer stays within that row. 

Eventually, since
the number of vertices in the bottom row is not a multiple of $3$, 
Minimizer is forced to use the cross-edge. At that point, both rows have 
number of vertices left a multiple of $3$; thus we are back in Case~1. 
\end{Proof}

\begin{theorem}
Consider the stripe-game, with Minimizer responding. Then  
 $R_3$ is perfect but $R_m$ for $m\ge 6$ is not.
\end{theorem}
\begin{Proof}
It can easily be checked that $R_3$ is perfect: whatever the first move, what is left is
connected on $3$ vertices and thus can be taken.

Consider $R_m$ for $m\ge 6$.
Minimizer's strategy will ensure that, until the very end, every vertex in the 
smaller row still has its cross-edge. Therefore, we can refer to the 
situation by just the counts of the two rows. We will use $(i,j)$, with $i\ge j$, to denote
the situation where one row has $i$ vertices and one row has $j$ vertices.
For the base of the induction, we need the case $(4,2)$. For this, one can readily
check that wherever Maximizer initiates, Minimizer can respond and 
disconnect the graph, thereby ending the game.

For $R_m$ 
the play starts at the case $(m,m)$. For the first move, Minimizer uses a 
cross-edge; so the case becomes $(m-1,m-2)$. 
We claim that Minimizer can ensure the case $(m-2,m-4)$ next. For, 
if Maximizer chooses a vertex in the larger side, then Minimizer stays in the larger side,
using up the vertex that has no cross-edge; and if 
Maximizer chooses a vertex in smaller side, then Minimizer takes two vertices there
and a cross-edge to the other row. 

We claim that Minimizer can ensure the case $(m-4,m-5)$ next. For, 
if Maximizer initiates on the larger side, then Minimizer takes three vertices there;
and if Maximizer initiates on the smaller side, then Minimizer immediately uses
the cross-edge and then takes one of the (two) vertices without a cross-edge.
By repeated application of the strategy, Minimizer can alternate between cases
of the form $(x,x-2)$ and $(y,y-1)$ until they reach the case $(4,2)$, which we saw is not perfect.
\end{Proof}

\section{The Unrooted $P_3$}

There is also a version of the game where the packing subgraph has no root.
We define the \dd{unrooted-$P_3$-game} to be the game where responder
need only choose a copy of $P_3$ containing the designated vertex. We show that
in some cases the value of the game is the same as in the rooted version, but in
other cases it is different. In particular, for Maximizer responding,
they get all that is possible.

\subsection{Maximizer responding}

\begin{theorem}
Consider the unrooted-$P_3$-game with Maximizer responding played on graph $G$.
Then the value of the game is $\mu(G)$.
\end{theorem}
\begin{Proof}
Consider a maximum $P_3$-packing $\curlyP$ of $G$. If Minimizer chooses
a vertex~$u$ in~$\curlyP$, then Maximizer responds with the associated
copy in $\curlyP$ and repeats. If instead Minimizer chooses a vertex $v$ outside $\curlyP$,
then by requirement, the vertex $v$ is in a copy $Q$ of~$P_3$. The only way 
a problem could arise is if $Q$ intersects two copies
 in~$\curlyP$. But that implies that vertex $v$ has an edge to some copy in $\curlyP$;
and so Maximizer can use that edge and one edge from that copy to
build a $P_3$, thereby affecting only one copy in~$\curlyP$. Repeat.
\end{Proof}

The above theorem generalizes Theorem~1 on the matcher game from~\cite{GH}.
Note that this pattern does not continue much further. In particular,
the star $K_{1,3}$ does not have a similar result. For example, take three copies
of~$K_{1,3}$ and add a new vertex $v$ adjacent to one end-vertex 
from each copy. The resulting tree has three disjoint $K_{1,3}$'s, but Minimizer
ends the game in one move by initiating on~$v$.

\subsection{Minimizer responding}

When one changes to the unrooted
game, this gives both player more options. We saw above that when Maximizer
is responder, the added options to Minimizer do not help them. A similar result
holds if we go from the stripe-game to the unrooted game with Minimizer as responder:

\begin{theorem}
Consider Minimizer responding. The value of the unrooted-$P_3$-game is 
at most the value of the stripe-game.
\end{theorem}
\begin{Proof}
Consider changing from the stripe-game to the unrooted game.
Minimizer plays just as if it were the stripe-game.
We argue that the new options do not help Maximizer. For, the only additional
option they have is to play a vertex $v$ that is in a $P_3$ but is not the end
of one. That  is, the only additional initiation option they get is
to play the center vertex of a star component; but that is equivalent
to initiating at a leaf of the component, which they could do already.
\end{Proof}

In contrast, the star-game is incomparable with the unrooted-$P_3$-game.
Consider, for example, the double corona of a complete graph (see Figure~\ref{f:doubleCorona}). 
Maximizer now can initiate on a leaf, and thereby ensure that
approximately half the vertices are used. On the other hand, consider the tree shown in Figure~\ref{f:nineTree}.
As we saw, Maximizer can obtain every vertex in the star-game by initiating on the central vertex $v$.
In the unrooted $P_3$-game, however, Minimizer can respond differently and destroy the perfection.

\section{Questions}

We conclude with some questions for future study. Obviously, a natural direction is to 
replace $P_3$ by another required subgraph. For the games with $P_3$, it would be
interesting to determine the value of the game on a general grid and a general rooks graph.
Another question, is whether there is a $\varepsilon > 0$ such that all graphs of order $n$ with minimum
degree at least $(1-\varepsilon)n$ are perfect.


\begin{thebibliography}{9}
\bibitem{GH}
W. Goddard and M.A. Henning.  The matcher game played in graphs.
 Discrete Appl. Math.  237  (2018), 82--88.
\bibitem{KKN}
A. Kaneko, A. Kelmans, and T. Nishimura. On packing $3$-vertex paths in a graph.
 J. Graph Theory  36  (2001),  no. 4, 175--197.
 \bibitem{KH} 
D.G. Kirkpatrick and P. Hell.  On the complexity of general graph factor problems.
 SIAM J. Comput.  12  (1983),  no. 3, 601--609.
 
\end{thebibliography}
\end{document}